\documentclass[11pt]{article}
\usepackage{amsmath}
\usepackage{graphicx}
\usepackage{epsfig}
\usepackage{amsfonts}
\usepackage{amssymb}
\usepackage{placeins}
\usepackage{cases}
\usepackage[margin=1in]{geometry}
\usepackage{authblk}
\usepackage[titletoc,toc,title]{appendix}
\usepackage{epstopdf}
\usepackage{changes}
\usepackage{xcolor}
\usepackage{bm}
\usepackage{enumerate}
\usepackage{caption}
\usepackage{subfig}
\usepackage{relsize}
\usepackage{amsmath}

\def\horizontaldistance{\kern2pt}
\def\verticaldistance{\kern 5pt}
\newcommand*{\vertbar}{\rule[-1ex]{0.5pt}{2.5ex}}

\usepackage{accents}

\usetikzlibrary{patterns}
\usepackage{tikz}
\usepackage[siunitx]{circuitikz}
\usetikzlibrary{decorations.pathreplacing,calligraphy}
\usetikzlibrary{patterns}
\usetikzlibrary{backgrounds}
\usetikzlibrary{decorations.text}
\usepackage{cite}
\usepackage{yfonts}
\usepackage{amsmath,amssymb,amsfonts}
\usepackage{tikz}

\usepackage{tikz}

\usetikzlibrary{patterns}
\usetikzlibrary{backgrounds}
\usetikzlibrary{decorations.text}

\usetikzlibrary{positioning,quotes,calc,fit,shapes,shadows,arrows,trees,mindmap}
\tikzset{block/.style={draw, very thick, minimum height=4cm, align=center}, line/.style={-latex}}
\tikzset{blockV/.style={draw, very thick, text width=2cm, minimum height=2cm, minimum width=4cm, align=center}, line/.style={-latex}}
\tikzset{blockExt/.style={draw, very thick, minimum height=0.7cm, minimum width=0.7cm, align=center}, line/.style={-latex}}

\usepgflibrary{arrows}
\usetikzlibrary{bayesnet}

\definecolor{color_gr}{RGB}{10, 120, 5}
\definecolor{color_gray}{rgb}{0, 0.05, 0.05}
\colorlet{color_vl}{violet!70}
\definecolor{light-gray}{HTML}{E0E0E0}
\definecolor{blue-violet}{rgb}{0.54, 0.17, 0.89}
\definecolor{light-gray}{HTML}{E0E0E0}
\definecolor{carnelian}{rgb}{0.7, 0.11, 0.11}
\definecolor{darkpastelgreen}{rgb}{0.01, 0.75, 0.24}

\title{Uncertainty Quantification in Data-Driven Dynamical Models via Inverse Problem Solving}

\begin{document}
\author[1]{Mohamed Akrout}
\author[1]{Dan Wilson}
\affil[1]{Department of Electrical Engineering and Computer Science, University of Tennessee, Knoxville, TN 37996, USA}
\maketitle

\begin{abstract}
Data-driven model identification strategies can be used to obtain phenomenological models that capture the temporal evolution of observable data.  While it is usually straightforward to obtain such a model from time series data, for instance with least-squares fitting, it is generally difficult to quantify the uncertainty associated with the prediction of the temporal evolution of the observables.  This paper considers a general framework for uncertainty quantification in data-driven dynamical models by framing prediction error through the lens of inverse problem theory.  Building on Koopman-inspired model identification strategies that are suited for nonlinear dynamical models,  we consider a prediction as an approximate measurement from which the original input state can be faithfully recovered, and define the prediction error as the MSE of solving the inverse problem that would yield this prediction. We demonstrate the efficacy of this approach on both numerical models and experimental data showing that it provides a robust uncertainty measure of model performance.
\end{abstract}

\section{Introduction}

Model identification is an imperative first step for analysis, control, and prediction of dynamical systems.  Data-driven model identification strategies, which yield phenomenological models that capture the temporal evolution of observable data (instead of physics-based models) can be particularly useful in many practical applications where the scale and complexity obscures information about the underlying physics and governing equations.  A wide variety of data-driven model identification strategies have been developed in recent years.  For instance, learning-based approaches can be used to identify sparse dynamical models that capture the evolution of observables \cite{mang19}, \cite{pant19}, \cite{brun16}.  In a similar manner, approaches using reservoir computing can be used to predict the temporal evolution of observables by using a fixed dynamical system (the reservoir) and learning an appropriate mapping between the states of the reservoir and the system observables \cite{gaut21}, \cite{tana19}.  Additionally, using Koopman operator theory as a foundation, dynamic mode decomposition can be used to obtain (generally) high-dimensional linear models that capture the evolution of observables on short-to-moderate time scales \cite{rowl09}, \cite{arba17}, \cite{schm10}, \cite{will15}.

While the aforementioned model identification approaches can be used to predict the temporal behavior of general nonlinear dynamical models, the presence of noise and other uncertainties can significantly degrade the predictive performance of the learned models, often necessitating an explicit strategy for handling these uncertainties.  One approach to working with noisy training data is to implement approaches to remove systematic biases due to noise \cite{hema17}, \cite{askh18}, \cite{than23}, \cite{went23}.  Alternatively, strategies have been developed to incorporate uncertainty and noise estimates into the model identification step, for instance, using bootstrap aggregating (bagging) methods \cite{fase22}, \cite{sash22} or by inferring both the underlying dynamical equations and measurement noise simultaneously \cite{rudy19}, \cite{kahe22}.  Bayesian approaches have also been applied to handle uncertainty in the model identification step \cite{nive24}, \cite{yang20}, \cite{green15}.  When applied to model identification of dynamical systems, these strategies can be used to obtain both a prediction of the future state and a prediction error that captures the uncertainty associated with that estimation.

When inferring dynamical models from data, prediction errors tend to accumulate over large prediction horizons representing a significant impediment to obtaining accurate models. While standard techniques focus on the estimation of the operator to evolve the state, they often lack a formal mechanism to quantify the reliability of the resulting predicted trajectory. We propose a framework that treats the prediction sequence as an observation in an inverse problem to derive uncertainty metrics.  In this work we consider the problem of state prediction for a general dynamical system from the perspective of inverse problem theory.  Leveraging prior Koopman-inspired frameworks for model identification for use in general nonlinear dynamical systems \cite{wils23nonlin}, this work casts dynamic prediction as an inverse problem and leverages Bayesian inference to obtain a measure of uncertainty based on the estimated second-order moment.   The proposed framework is built on the principle of self-consistency:~if a data-driven predictive model is reliable, it should remain identifiable when subjected to an inverse mapping. By treating the predicted states as pseudo-measurements, we can evaluate how much information is preserved or lost during the forward evolution.  We treat the model’s own forward prediction as a pseudo-measurement and subsequently invert it. Here, the prediction serves as the ground truth for the inverse problem which aligns with several established paradigms in signal processing and machine learning that utilize self-referential feedback to estimate reliability \cite{sutton1988learning,hinton2006reducing,zhu2017unpaired}.

The organization of this paper is as follows:~Section \ref{backsec} provides relevant background on recently developed Koopman-based model identification algorithms as well as the vector approximate message passing (VAMP) algorithm that is used in this work to solve linear inverse problems.   Section \ref{methodsec} describes the proposed algorithm describing the approach of quantifying error prediction through Bayesian inverse reformulation.  Section \ref{ressec} demonstrates this approach on both numerical models and experimental data illustrating its utility.  Section \ref{concsec} gives concluding remarks.

\section{Background} \label{backsec}

\subsection{Koopman Operator Theory}
Consider a discrete-time dynamical system of the form
\begin{equation} \label{disceqn}
    \bm{x}^{k+1} = F(\bm{x}^k),
\end{equation}
where $\bm{x}^k \in \mathbb{R}^n$ is the state and $F$ determines the dynamics of the mapping $\bm{x}^k \mapsto \bm{x}^{k+1}$.  Let $\psi:\mathbb{R}^n \rightarrow \mathbb{F}$ be an observable where $\psi \in \mathcal{F}$ is the observable space.  The Koopman operator $K:\mathcal{F} \rightarrow \mathcal{F}$ is defined as
\begin{equation} \label{kaut}
    K \psi(\bm{x}^k) \equiv \psi(F(\bm{x}^k)),
\end{equation}
and can be thought of as giving the dynamics of observables for the system \eqref{disceqn}.  The Koopman operator is linear, inheriting this property from the linearity of the composition \cite{budi12}, \cite{mezi13}.  Time-varying external inputs can also be considered using Koopman-based approaches \cite{kord18}, \cite{will16}, \cite{schm10}, starting with a system of the form
\begin{equation} \
     \bm{x}^{k+1} = F(\bm{x}^k,\bm{u}^k),
\end{equation}
where $\bm{u}^k \in \mathcal{U} \subset \mathbb{R}^m$.  As in \cite{kord18}, letting $l(\mathcal{U}) = \{ (\bm{u}_k)_{k = 0}^\infty | \bm{u}_k \in \mathcal{U} \}$ be the space of all input sequences, and defining an observable $\psi:\mathbb{R}^n \times l(\mathcal{U}) \rightarrow \mathbb{R}$, a non-autonomous version of the Koopman operator is
\begin{equation} \label{knaut}
    K \psi(\bm{x}^k,(\bm{u}^k)_{k = 0}^\infty)) = \psi(F(\bm{x}^k,\bm{u}^0), (\bm{u}^k)_{k = 1}^\infty ).
\end{equation}
In both the autonomous \label{kaut} and non-autonomous \label{knaut} case, the Koopman operator is linear but generally infinite dimensional. 

\subsection{Data-Driven Model Identification Using Koopman Operator Theory} \label{kooptheory}

The primary challenge to using this approach in data-driven systems is in finding a finite-dimensional approximation to this linear operator.  Dynamic mode decomposition is one such approach \cite{kutz16}, \cite{schm10}, \cite{rowl09}, whereby a series of snapshots is collected
\begin{equation}
    s_k = (g(\bm{x}^k),g(\bm{x}^{k+1})),
\end{equation}
defined for $k = 1,\dots,q$ where $g(\bm{x}^k) \in \mathbb{R}^p = \begin{bmatrix} \psi_1,\dots,\psi_p \end{bmatrix}$ is a set of $p$ observables.  Arranging the snapshot pairs into matrices $\bm{X} = \begin{bmatrix} g(\bm{x}^1) & \dots & g(\bm{x}^{q-1})  \end{bmatrix}$ and $\bm{X}^+ = \begin{bmatrix} g(\bm{x}^2) & \dots & g(\bm{x}^{q})  \end{bmatrix}$, one can obtain a finite dimensional approximation of the Koopman operator according to 
\begin{equation}
    g(\bm{x}^{k+1}) = \bm{A} g(\bm{x}^k),
\end{equation}
where $\bm{A} = \bm{X}^+ \bm{X}^\dagger$ and $^\dagger$ denotes the pseudoinverse.   Many extensions of DMD have been proposed including extended DMD (EDMD) \cite{will15} which considers a lifted observable space, and DMD with control (DMDc) \cite{kord18}, \cite{proc16} which explicitly considers a non-autonomous input in the model identification, and Hankel DMD \cite{arba17} which considers a delay embedding as an initial lifting of the observable space.  As an example of such an approach, drawing from the Hankel DMD approach \cite{proc16}, one can consider a time-delay embedding of the input and observables
\begin{equation}\label{deq}
    \bm{h}^k = \begin{bmatrix}  g(\bm{x}^{k-1})  \\ \vdots \\ g(\bm{x}^{k-z}) \\ \bm{u}^{k-1} \\ \vdots \\\bm{u}^{k-z}  \end{bmatrix}.
\end{equation}
Here $\bm{h}^k \in \mathbb{R}^M$ where $M = z(p+m)$.  Drawing from the extended DMD approach \cite{will15}, the delay embedding from \eqref{deq} can be lifted to a higher dimensional space according to
\begin{equation} \label{liftfn}
    \bm{\upsilon}^k = f_{\rm lift}(g(\bm{x}^k),\bm{h}^k),
\end{equation}
where $f_{\rm lift}(g(\bm{x}^k),\bm{h}_k) \in \mathbb{R}^L$ with $L > M$ is a general nonlinear function.  Note that from the perspective of Koopman operator theory, $\begin{bmatrix} g(\bm{x}^k)^T  & (\bm{\upsilon}^k)^T   \end{bmatrix}^T$ can be viewed as an observable for the system, and letting $\bm{\Upsilon} = \begin{bmatrix} \bm{\upsilon}^1 & \dots &  \bm{\upsilon}^{q-1} \end{bmatrix}$, $\bm{\Upsilon}^+ = \begin{bmatrix} \bm{\upsilon}^2 & \dots & \bm{\upsilon}^{q} \end{bmatrix}$, and $\bm{U} = \begin{bmatrix} \bm{u}^1 & \dots & \bm{u}^{q-1} \end{bmatrix}$, a finite dimensional approximation of the associated Koopman operator, $\bm{A}_\Upsilon$ follows
\begin{equation}
    \begin{bmatrix} \bm{X}^+ \\ \bm{\Upsilon}^+  \end{bmatrix} = \bm{A}_\Upsilon   \begin{bmatrix} \bm{X} \\ \bm{U} \\ \bm{\Upsilon}  \end{bmatrix},
\end{equation}
which can be approximated according to 
\begin{equation}
    \bm{A}_\Upsilon =  \begin{bmatrix} \bm{X}^+ \\ \bm{\Upsilon}^+  \end{bmatrix}   \begin{bmatrix} \bm{X} \\ \bm{U} \\ \bm{\Upsilon}  \end{bmatrix}^\dagger,
\end{equation}
where $\bm{A}_\Upsilon \in \mathbb{R}^{(p+L) \times (p+L+m) }$.  Subsequently, 
\begin{equation} \label{linkoop}
    \begin{bmatrix}  g(\bm{x}^{k+1}) \\ \bm{\upsilon}^{k+1}   \end{bmatrix} = \bm{A}_\Upsilon   \begin{bmatrix}  g(\bm{x}^k) \\  \bm{u}^k \\ \bm{\upsilon}^{k} \end{bmatrix} 
\end{equation}
can be used to predict the evolution of the observables $g(\bm{x}^k)$ and the lifted state.  This is similar to the approach used in \cite{kord18}.

\subsection{Nonlinear Estimators for the Koopman Operator} \label{nonlinest}
In the previous section \ref{linkoop} provides a linear estimator for the action of the Koopman operator on the lifted observable space.  Note, however, that it is generally not necessary predict the evolution of the variables $\bm{\upsilon}^k$ since they can be uniquely determined directly from measurements of the observables.  The approach presented in \cite{wils23nonlin} considers a related method for finding nonlinear estimators for the action of the Koopman operator, which can ultimately provide low-dimension, data-driven models that accurately reflect the behavior of highly nonlinear systems. In this case, a data-driven model is sought that follows
\begin{align} \label{nlinest}
    g(\bm{x}^{k+1}) &= \bm{A}_N  \begin{bmatrix}  g(\bm{x}^k) \\  \bm{u}^k \\ \bm{\upsilon}^{k} \end{bmatrix},   \nonumber \\
    \bm{h}^{k+1} &= f_h(\bm{h}^k,g(\bm{x}^{k+1}),\bm{u}^k), \nonumber \\
    \bm{\upsilon}^{k+1} &= f_{\rm lift}(g(\bm{x}^k),\bm{h}^{k+1}),
\end{align}
where $\bm{A}_N \in \mathbb{R}^{p \times (p+L+m)}$, $f_h$ shifts the entries from $\bm{h}^k$ defined in \eqref{deq} with the incoming data, and $f_{\rm lift}$ is the same nonlinear function defined in \eqref{liftfn}.  The matrix $\bm{A}_N$ is the only unknown in the above equation, and a least-squares estimate can be obtained according to
\begin{equation}
     \bm{A}_N =   \bm{X}^+     \begin{bmatrix} \bm{X} \\ \bm{U} \\ \bm{\Upsilon}  \end{bmatrix}^\dagger,
\end{equation}
The size of $\bm{A}_N$ is substantially smaller than $\bm{A}_\Upsilon$ and hence, requires the inference of far fewer coefficients, mitigating the risk of overfitting.  Of course, Equation \eqref{nlinest} is nonlinear and compared to \eqref{linkoop} cannot be analyzed simply in terms of an eigendecomposition.  Nonetheless, as seen in \cite{wils23nonlin}, Equation \eqref{nlinest} can be used to replicate fundamentally nonlinear behaviors (e.g.,~limit cycles, bistability), and generally provide more accurate predictions on substantially longer timescales than the corresponding linear models \eqref{linkoop}.

The number of coefficients required in the fitting of \eqref{nlinest} can further be reduced by considering the data used for fitting in terms of a reduced order basis.  Proper orthogonal decomposition (POD) \cite{holm96} is one common strategy used for this purpose.  Letting $\bm{B} = \begin{bmatrix} \bm{X}^T & \bm{U}^T & \bm{\Upsilon}^T \end{bmatrix}^T$, POD modes can be found according to the eigenvectors, $\bm{v}_1,\bm{v}_2,\dots$ of $\bm{B} \bm{B}^T$ and sorted in descending order according to the magnitude of their associated eigenvalues.  The first $\zeta$ POD modes associated with the largest eigenvalues are used to define an orthogonal basis $\bm{\Phi} = \begin{bmatrix} \bm{v}_1 & \dots & \bm{v}_\zeta   \end{bmatrix} \in \mathbb{R}^{ (p + L + m) \times \zeta}$ for which
\begin{equation}
    \begin{bmatrix}  g(\bm{x}^k) \\  \bm{u}^k \\ \bm{\upsilon}^{k} \end{bmatrix} \approx \sum_{j = 1}^\zeta  \bm{v}_j \alpha_{j,k},
\end{equation}
where $\alpha_{j,k} =  \begin{bmatrix} g(\bm{x}^k)^T  & (\bm{u}^k)^T & (\bm{\upsilon}^k)^T   \end{bmatrix} \bm{v}_j$  gives the contribution the $j^{\rm th}$ POD mode.  The model fitting can then be done instead in the POD basis, in this case the first line in \eqref{nlinest} can be written in a 
\begin{equation}\label{eq:predicition-formula}
   g(\bm{x}^{k+1}) = \bm{A}_R  \bm{\Phi}^T  \begin{bmatrix}  g(\bm{x}^k) \\  \bm{u}^k \\ \bm{\upsilon}^{k} \end{bmatrix}.
   \end{equation}
and a least-squares estimate for $\bm{A}_R \in \mathbb{R}^{p \times \zeta}$ can be obtained according to 
\begin{equation}
    \bm{A}_R = \bm{X}^+ \left(  \bm{\Phi}^T  \begin{bmatrix} \bm{X} \\ \bm{U} \\ \bm{\Upsilon}  \end{bmatrix}  \right) ^\dagger.
\end{equation}

\subsection{Inverse problems and expectation propagation}\label{sec:EP}

Inverse problems arise when one seeks to recover an unknown signal 
$\bm{x} \in \mathbb{R}^N$ from indirect and noisy observations of the form
\begin{equation}\label{eq:linear-model}
    \bm{y} = \bm{A} \bm{x} + \bm{n},
\end{equation}
where $\bm{A} \in \mathbb{R}^{M \times N}$ is a known forward operator and $\bm{n} \in \mathbb{R}^{M}$ represents measurement noise. This formulation appears in numerous scientific and engineering domains such as compressed sensing, tomography, imaging, and wireless communications. The central challenge lies in the fact that the forward operator $\bm{A}$ may be ill-conditioned or rank-deficient, making the inversion problem inherently unstable in the presence of noise.  

From a Bayesian perspective, inverse problems can be framed as estimating 
the posterior distribution $p(\bm{x}|\bm{y})$ which is given by the Bayes rule as
\begin{equation}
    p(\bm{x}|\bm{y}) \propto p(\bm{y}|\bm{x})\,p(\bm{x}),
\end{equation}
where $p(\bm{y}|\bm{x})$ encodes the likelihood model (often Gaussian, corresponding  to additive white noise) and $p(\bm{x})$ represents prior knowledge about the signal structure (e.g., sparsity, smoothness, or non-negativity). This probabilistic formulation turns the recovery task into Bayesian inference, where one seeks point estimates such as the minimum mean square error (MMSE) or the maximum a posteriori (MAP) or posterior mean.

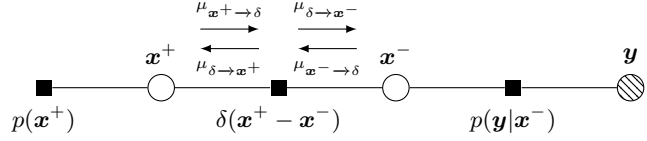
\begin{figure}[!h]
    \centering
    \begin{tikzpicture}[scale = 1.3]
    \newcommand{\vertex}{\node[vertex]}
    \tikzset{vertex/.style = {circle, draw, inner sep = 0pt, minimum size = 10pt}}
    \vertex[label = $\bm{x}^+$](x1) at (-1.2,0) {};
    \vertex[label = $\bm{x}^-$](x2) at (1.2,0) {};
    \vertex[label = {$\bm{y}$}, pattern = north west lines](y) at (3.6,0) {};
    \tikzset{vertex/.style = {rectangle, fill = black, inner sep = 0pt, minimum size = 6pt}}
    \vertex[label = below: $p(\bm{x}^+)$](px) at (-2.4,0) {};
    \vertex[label = below: $\delta(\bm{x}^+ - \bm{x}^-)$](delta) at (0,0) {};
    \vertex[label = below: {$p(\bm{y}|\bm{x}^-)$}](lik) at (2.4,0) {};
    \draw (px)--(x1);
    \draw (delta)--(x1);
    \draw (delta)--(x2);
    \draw (lik)--(x2);
    \draw (lik)--(y);
    \draw[-latex] ([xshift=0.2cm,yshift=0.6cm] delta.center)--([xshift=-0.4cm,yshift=0.6cm] x2.center);
    \node[] at ([xshift=1.7cm,yshift=0.8cm] x1.center) {\tiny{$\mu_{\delta \xrightarrow \,\bm{x}^-}$}};
    \draw[-latex] ([xshift=-0.4cm,yshift=0.4cm] x2.center)--([xshift=0.2cm,yshift=0.4cm] delta.center);
    \node[] at ([xshift=1.7cm,yshift=0.2cm] x1.center) {\tiny{$\mu_{\bm{x}^- \xrightarrow \,\delta}$}};
    \draw[-latex] ([xshift=-0.8cm,yshift=0.6cm] delta.center)--([xshift=-1.4cm,yshift=0.6cm] x2.center);
    \node[] at ([xshift=0.7cm,yshift=0.8cm] x1.center) {\tiny{$\mu_{\bm{x}^+ \xrightarrow \,\delta}$}};
    \draw[-latex] ([xshift=-1.4cm,yshift=0.4cm] x2.center)--([xshift=-0.8cm,yshift=0.4cm] delta.center);
    \node[] at ([xshift=0.7cm,yshift=0.2cm] x1.center) {\tiny{$\mu_{\delta \xrightarrow \,\bm{x}^+}$}};
\end{tikzpicture}
    \caption{VAMP's factor graph according to the factorization (\ref{eq:VAMP_posterior-factorize}).}
    \label{fig:VAMP-factor-graph-original}
\end{figure}

The vector approximate message passing (VAMP) algorithm is designed to solve linear inverse problem in (\ref{eq:linear-model}) under known Gaussian or non-Gaussian priors $p(\bm{x})$, overcoming the computational intractability of MMSE estimation. VAMP resorts to an approximation of the belief propagation algorithm, namely expectation propagation (EP) \cite{minka2001family}. EP approximates the true posterior distribution of $\bm{x}$ by iteratively matching moments (mean and variance) between the true posterior and an exponential family distribution, typically Gaussian. This moment-matching process ensures that the approximation captures key global characteristics of the posterior, rather than being limited to local optima or mode-centric approximations. Unlike methods requiring strict Gaussian assumptions, EP can handle non-Gaussian priors on evolving states and likelihoods by decomposing the joint distribution through local factors forming the so-called ``factor graph'' shown in Fig.~\ref{fig:VAMP-factor-graph-original}. Specifically, after splitting $\bm{x}$ into two separate variables $\bm{x}^+$ and $\bm{x}^-$, the conditional joint density of $\bm{x}^+$ and $\bm{x}^-$ given $\bm{y}$ is factorized as follows:
\begin{equation}
\label{eq:VAMP_posterior-factorize}
p(\bm{x}^+, \bm{x}^- | \bm{y})\propto p(\bm{x}^+)\, \delta(\bm{x}^+ - \bm{x}^-)\, p(\bm{y}|\bm{x}^-).
\end{equation}

Each non-Gaussian factor is approximated locally by a Gaussian term while preserving the overall non-Gaussian structure, making it suitable for complex dynamical systems with nonlinear state transitions and non-Gaussian noise. Furthermore, by constraining the approximate posterior to a Gaussian form, EP yields a tractable distribution where mean and variance can be computed analytically. This enables efficient inference, prediction, and uncertainty quantification without requiring computationally intensive sampling methods.

\begin{figure*}[t!]
        \centering
         \begin{tikzpicture}[thick,scale=0.7, every node/.style={transform shape}]
  \node[block, fill=blue!15, minimum width=2.5cm,text width=2.5cm] (z_uv) {\Huge{$p(\boldsymbol{x}^+)$}};
  \node[block, fill=green!15, minimum width=5cm,text width=5cm,right= 5cm of z_uv] (phi_z) {\huge{
  \mbox{$\boldsymbol{y} = \bm{A}\,\boldsymbol{x}^- + \bm{n}$}}};
  \node[blockExt,right=of z_uv, xshift=0.4cm, yshift=-1.7cm] (ext_z_uv_to_y) {$\mathrm{\textbf{EXT}}$};
  \node[blockExt,left=of phi_z, xshift=-0.4cm, yshift=1.6cm] (ext_y_to_z_uv) {$\mathrm{\textbf{EXT}}$};
  
  \draw [-latex,very thick] ([yshift=-4.4em]z_uv.east) --
  node [midway,below=0em,align=center ] { \Large{$\widehat{\boldsymbol{x}}_{\mathsf{p}}^{+}$}}
  node [midway,below=1.6em,align=center ] {\Large{$\gamma_{\boldsymbol{x}_{\mathsf{p}}^{+}}$}}
  (ext_z_uv_to_y.west);
  \draw [-latex,very thick] (ext_z_uv_to_y) --
  node [midway,below=0.1em,align=center ] { \Large{$\widehat{\bm{x}}^-_{\mathsf{e}}$}}
  node [midway,below=1.7em,align=center ] {\Large{$\gamma_{\bm{x}_{\mathsf{e}}^-}$}}
  ([yshift=-4.4em]phi_z.west)
  node [pos=0.25](ext_between_z_uv_y){};
  \draw [-latex,very thick] ([yshift=4.2em]phi_z.west) --
  node [midway,above=0em,align=center ] { \Large{$\widehat{\boldsymbol{x}}_{\mathsf{p}}^{-}$}}
  node [midway,above=1.8em,align=center ] {\Large{$\gamma_{\boldsymbol{x}_{\mathsf{p}}^{-}}$}}
  (ext_y_to_z_uv.east);
  \draw [-latex,very thick] (ext_y_to_z_uv) --
  node [midway,above=0em,align=center ] { \Large{$\widehat{\boldsymbol{x}}_{\mathsf{e}}^{+}$}}
  node [midway,above=1.7em,align=center ] {\Large{$\gamma_{\boldsymbol{x}_{\mathsf{e}}^+}$}}
  ([yshift=4.25em]z_uv.east)
  node [pos=0.25](ext_between_y_z_uv){};
  \draw [-latex,very thick] ([xshift=-0.5em]ext_between_z_uv_y.center) --
  (ext_y_to_z_uv.south);
  \draw [-latex,very thick] ([xshift=0.5em]ext_between_y_z_uv.center) --
  (ext_z_uv_to_y.north);

\end{tikzpicture}
         \vspace{-0.1cm}
        \caption{Block diagram of the VAMP approach to solve the inverse problem with its two modules: one MMSE denoiser incorporating the prior $p(\bm{x}^+)$, and the LMMSE module. The two modules exchange extrinsic information/messages through the \protect\tikz[inner sep=.25ex,baseline=-.75ex] \protect\node[rectangle,draw,thick,minimum width=0.45cm,minimum height=0.45cm] {\footnotesize \textbf{EXT}}; blocks.}
        \label{fig:block-diagram}
\end{figure*}
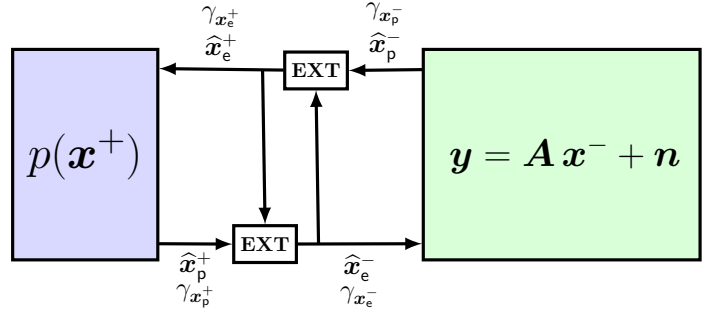

\noindent The VAMP algorithm runs by iterating between two estimation steps depicted in Fig. \ref{fig:block-diagram}: $i)$ enforcing the prior $p(\bm{x})$ through MMSE estimation, and $ii)$ enforcing the channel model $p(\bm{y}|\bm{x})$ given the observed measurement $\bm{y}$. In the sequel, we delve into the derivation details of each step.

\paragraph{Step 1: MMSE estimation of $\bm{x}^+$}

At iteration $k$, a separable prior $p(\bm{x}^+) = \prod_{n=1}^{N} p(x_n^+)$ and the incoming extrinsic Gaussian belief from the node $\delta(\bm{x}^+-\bm{x}^-)$ are available. The latter takes the form $\mu_{\delta\to\bm{x}^+} = \mathcal{CN}\big(\bm{x}^+; \widehat{\bm{x}}_{\textsf{e},k-1}^+, \gamma_{\bm{x}_{\textsf{e}}^+,k-1}^{-1}\mathbf{I}_{N}\big)$, which can be interpreted as a pseudo‑observation, $\widehat{\bm{x}}_{\textsf{e},k-1}^+$, of the underlying signal $\bm{x}^+$, corrupted by i.i.d. AWGN with variance $\gamma_{\bm{x}_{\textsf{e}}^+,k-1}^{-1}$. Consequently, the corresponding likelihood is $p(\widehat{\bm{x}}_{\textsf{e},k-1}^+|\bm{x}^+) = \mu_{\delta\to\bm{x}^+}$. Using EP, the posterior belief of $\bm{x}^+$ is approximated as Gaussian: $\mathcal{CN}\big(\bm{x}^+; \widehat{\bm{x}}_{\textsf{p},k}^+, \gamma_{\bm{x}_{\textsf{p}}^+,k}^{-1} \mathbf{I}_N\big)$. The posterior mean $\widehat{\bm{x}}_{\textsf{p},k}^+$ and variance $\gamma_{\bm{x}_{\textsf{p}}^+,k}^{-1}$ are given by \cite{rangan2019vector}:
\begin{subequations}
\label{eq:VAMP_pst-belief-x+}
\begin{align}
    \widehat{\bm{x}}_{\textsf{p},k}^+ &= \int \bm{x}\, p\big(\bm{x}|\widehat{\bm{x}}_{\textsf{e},k-1}^+\big)\, \textrm{d}\bm{x}\triangleq \bm{g}_{\bm{\mathsf{x}}^+}\big(\widehat{\bm{x}}_{\textsf{e},k-1}^+, \gamma_{\bm{x}_{\mathsf{e}}^+,k-1}\big),\label{eq:VAMP_pst-belief-x+_mean}\\
    \gamma_{\bm{x}_{\textsf{p}}^+,k}^{-1} &= \gamma_{\bm{x}_{\mathsf{e},k-1}^+}^{-1} \left\langle \bm{g}^{\prime}_{\bm{\mathsf{x}}^+}\big(\widehat{\bm{x}}_{\textsf{e},k-1}^+,\gamma_{\bm{x}_{\mathsf{e}}^+,k-1}\big)\right\rangle,\label{eq:VAMP_pst-belief-x+_precision}
\end{align}
\end{subequations}
    
\noindent with
\begin{equation}
    p_{\bm{\mathsf{x}}^+ \mid \bm{\mathsf{x}}_{\textsf{e},k-1}^+}\big(\bm{x}|\widehat{\bm{x}}_{\textsf{e},k-1}^+\big) \propto p_{\bm{\mathsf{x}}^+}(\bm{x})\, \mathcal{CN}\big(\bm{x}; \widehat{\bm{x}}_{\mathsf{e},k-1}^+ ,\gamma^{-1}_{\bm{x}_{\mathsf{e}}^+,k-1} \mathbf{I}_{N}\big),
\end{equation}

\noindent and
\begin{equation}
\label{eq:VAMP_denoiser-prime-x+}
    \bm{g}^{\prime}_{\bm{\mathsf{x}}^+}\big( \widehat{\bm{x}}_{\textsf{e},k-1}^+,\gamma_{\bm{x}_{\mathsf{e}}^+,k-1} \big) = \mathrm{diag}\bigg(\frac{ \partial\bm{g}_{\bm{\mathsf{x}}^+}\big(\widehat{\bm{x}}_{\textsf{e},k-1}^+,\gamma_{\bm{x}_{\textsf{e}}^+,k-1}\big) }{ \partial\widehat{\bm{x}}_{\textsf{e},k-1}^+ }\bigg).
\end{equation}

\noindent Note how the diagonal elements of the Jacobian in (\ref{eq:VAMP_denoiser-prime-x+}) are being averaged in (\ref{eq:VAMP_pst-belief-x+_precision}). This stems from approximating the covariance matrix with a scaled identity matrix. The EP belief approximation in (\ref{eq:VAMP_pst-belief-x+}) is a special variant of the EP algorithm \cite{minka2001family}, namely scalar EP \cite{cakmak2018expectation}. Thereafter, the outgoing extrinsic Gaussian belief $\mu_{\bm{x}^+\to\delta}$ is updated as $\mathcal{CN}\big(\bm{x}^+; \widehat{\bm{x}}_{\textsf{p},k}^+, \gamma_{\bm{x}_{\textsf{p}}^+,k}^{-1}\mathbf{I}_N\big) / \mathcal{CN}\big(\bm{x}^+; \widehat{\bm{x}}_{\textsf{e},k-1}^+, \gamma_{\bm{x}_{\textsf{e}}^+,k-1}^{-1}\mathbf{I}_{N}\big)$, and then sent back in the form of a mean vector and a scalar precision given by \cite{rangan2019vector}:
\begin{subequations}
\label{eq:VAMP_ext-belief-x-}
\begin{align}
    \gamma_{\bm{x}_{\textsf{e}}^-,k} &= \gamma_{\bm{x}_{\textsf{p}}^+,k} - \gamma_{\bm{x}_{\textsf{e}}^+,k-1}, \label{eq:VAMP_ext-belief-x-_mean}\\
    \widehat{\bm{x}}_{\textsf{e},k}^- &= \gamma_{\bm{x}_{\textsf{e}}^-,k}^{-1} \Big(\gamma_{\bm{x}_{\textsf{p}}^+,k}\, \widehat{\bm{x}}_{\textsf{p},k}^+ - \gamma_{\bm{x}_{\textsf{e}}^+,k-1}\, \widehat{\bm{x}}_{\textsf{e},k-1}^+\Big).\label{eq:VAMP_ext-belief-x-_precision}
\end{align}
\end{subequations}

\paragraph{Step 2: Linear MMSE (LMMSE) estimation of $\bm{x}^-$} 

Given the extrinsic Gaussian belief $\mu_{\delta\to\bm{x}^-} = \mathcal{CN}\big(\bm{x}^-; \widehat{\bm{x}}_{\textsf{e},k}^-, \gamma_{\bm{x}_{\textsf{e}}^-,k}^{-1}\mathbf{I}_N\big)$ and the measurement vector $\bm{y}$ according to the model (\ref{eq:linear-model}), the LMMSE estimator, $\bm{g}_{\bm{\mathsf{x}}^-}(\cdot)$, is obtained in closed form, along with its posterior variance as follows:
\begin{subequations}
\label{eq:VAMP_pst-belief-x-}
\begin{align}
    \widehat{\bm{x}}_{\textsf{p},k}^- &= \big( \gamma_w\, \bm{A}^{\textsf{H}}\bm{A} + \gamma_{\bm{x}_{\textsf{e}}^-,k}\, \mathbf{I}_N \big)^{-1} \big( \gamma_w\, \bm{A}^{\textsf{H}}\, \bm{y} + \gamma_{\bm{x}_{\textsf{e}}^-,k}\, \widehat{\bm{x}}_{\textsf{e},k}^- \big)\nonumber\\
    &\triangleq \bm{g}_{\bm{\mathsf{x}}^-}\big(\widehat{\bm{x}}_{\textsf{e},k}^-, \gamma_{\bm{x}_{\textsf{e}}^-,k}\big),\label{eq:VAMP_pst-belief-x-_mean}\\
    \gamma_{\bm{x}_{\textsf{p}}^-,k}^{-1} &= \tfrac{1}{N}\, \mathrm{Tr}\Big[ \big( \gamma_w\, \bm{A}^{\textsf{H}}\bm{A} + \gamma_{\bm{x}_{\textsf{e}}^-,k}\, \mathbf{I}_{N} \big)^{-1} \Big]\nonumber\\
    & \triangleq \gamma_{\bm{x}_{\textsf{e}}^-,k}^{-1}\, \Big\langle \bm{g}_{\bm{\mathsf{x}}^-}^{\prime}\big(\widehat{\bm{x}}_{\textsf{e},k}^-, \gamma_{\bm{x}_{\textsf{e}}^-,k}\big) \Big\rangle.\label{eq:VAMP_pst-belief-x-_gamma}
    \end{align}
\end{subequations}

\noindent with
\begin{equation}
    \bm{g}_{\bm{\mathsf{x}}^-}^{\prime}\big(\widehat{\bm{x}}_{\textsf{e},k}^-, \gamma_{\bm{x}_{\textsf{e}}^-,k}\big) = \mathrm{diag}\bigg(\frac{ \partial\bm{g}_{\bm{\mathsf{x}}^-}\big(\widehat{\bm{x}}_{\textsf{e},k}^-, \gamma_{\bm{x}_{\textsf{e}}^-,k}\big) }{ \partial\widehat{\bm{x}}_{\textsf{e},k}^- }\bigg).
\end{equation}

\noindent To mitigate the high computational cost stemming from the matrix inversion in (\ref{eq:VAMP_pst-belief-x-}), the singular value decomposition (SVD) of the matrix $\bm{A} = \bm{U}\, \mathrm{Diag}(\bm{s})\, \bm{V}^{\textsf{H}}$ is used to reduce the complexity of each VAMP iteration from $\mathcal{O}(N^3)$ to $\mathcal{O}(NR)$, with $R$ being the rank of $\bm{A}$. We refer the reader to \cite{rangan2019vector} for further details about the derivation of this step. The update of the outgoing extrinsic Gaussian belief $\mu_{\bm{x}^-\to\delta}$, similarly to the MMSE module, is accomplished through the computation of the extrinsic mean vector and scalar precision:
\begin{subequations}
\label{eq:VAMP_ext-belief-x+}
\begin{align}
    \gamma_{\bm{x}_{\textsf{e}}^+,k} &= \gamma_{\bm{x}_{\textsf{p}}^-,k} - \gamma_{\bm{x}_{\textsf{e}}^-,k}, \label{eq:VAMP_ext-belief-x+_mean}\\
    \widehat{\bm{x}}_{\textsf{e},k}^+ &= \gamma_{\bm{x}_{\textsf{e}}^+,k}^{-1} \Big(\gamma_{\bm{x}_{\textsf{p}}^-,k}\, \widehat{\bm{x}}_{\textsf{p},k}^- - \gamma_{\bm{x}_{\textsf{e}}^-,k}\, \widehat{\bm{x}}_{\textsf{e},k}^-\Big),\label{eq:VAMP_ext-belief-x+_precision}
\end{align}
\end{subequations}

\noindent which are sent back to the MMSE module for the next iteration.

In summary, solving inverse problems in (\ref{eq:linear-model}) are naturally 
addressed using Bayesian inference, especially using message passing algorithms that rely in the EP approximation. The latter is not only critical to ensure computational efficiency, but also to provide a second-moment information which will serve as the foundation for our proposed method for error quantification.

\section{Error Uncertainty Estimation via Inverse Problems}  \label{methodsec}

Existing approaches treat the predicted trajectory as a deterministic truth by neglecting the uncertainty inherent in the operator approximation. In this section, we describe a new estimation framework that recasts dynamic prediction as an inverse problem, and leverages Bayesian inference to obtain a measure of uncertainty based on the estimated second-order moment.

\subsection{Methodology}
\begin{figure}[h!]
    \centering
    \includegraphics[scale=0.62]{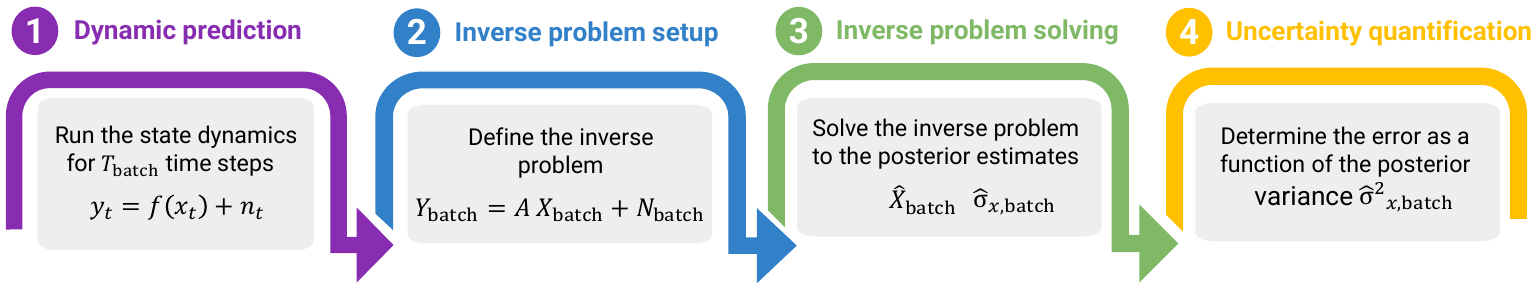}\vspace{-0.2cm}
    \caption{Quantification of prediction error through Bayesian inverse reformulation and solving}
    \label{fig:inverse-problem-framework}
\end{figure}

The key idea of the proposed framework is built on the principle of self-consistency: If a data-driven predictive model is reliable, it should remain identifiable when subjected to an inverse mapping. By treating the predicted states as ``pseudo-measurements'', we can evaluate how much information is preserved or lost during the forward evolution. The operational workflow of this framework is illustrated in Fig. \ref{fig:inverse-problem-framework}.

The process of error quantification begins with \textit{dynamic prediction}, where the estimated state dynamics are run for a window of $T_{\text{batch}}$ time steps. Following Eq. \ref{eq:predicition-formula}, the system generates a sequence estimated states $g(\bm{x}^{k+1})$ using the estimated data-driven model. Rather than accepting $g(\bm{x}^{k+1})$ at face value, we rewrite the predicted states obtained over $T_{\text{batch}}$ time steps as follows:
\begin{equation}\label{eq:Y=AX} 
\underbrace{\left[
  \begin{array}{cccc}
    \vertbar & \vertbar &        & \vertbar \\
    g(\bm{x}^{t_0+1})    & g(\bm{x}^{t_0+2})    & \ldots & g(\bm{x}^{T_{\textrm{batch}}+1})    \\
    \vertbar & \vertbar &        & \vertbar 
  \end{array}
\right]}_{\triangleq ~\bm{Y}_{\textrm{batch}}} = \underbrace{\bm{A}_R  \,\boldsymbol{\Phi}^\top}_{\triangleq ~\bm{A}} \underbrace{\left[
  \begin{array}{cccc}
    \vertbar & \vertbar &        & \vertbar \\
    g(\bm{x}^{t_0})    & g(\bm{x}^{t_0+1})    &  & g(\bm{x}^{T_{\textrm{batch}}})    \\
     \bm{u}^{t_0} & \bm{u}^{t_0+1}& \ldots & \bm{u}^{T_{\textrm{batch}}}\\
     \bm{\upsilon}^{t_0} & \bm{\upsilon}^{t_0+1}&  & \bm{\upsilon}^{T_{\textrm{batch}}}\\
    \vertbar & \vertbar &        & \vertbar 
  \end{array}
\right]}_{\triangleq ~\bm{X}_{\textrm{batch}}}.
\end{equation}

\noindent In this second step \textit{inverse problem setup}, we treat the predicted sequence of states $\bm{Y}_{\textrm{batch}}$ in (\ref{eq:Y=AX}) as an observation from which we attempt to recover the underlying state sequence $\bm{X}_{\text{batch}}$. By doing so, this formulation allows us to test the reliability of the prediction; an inaccurate data-driven model model will result in an observation $\bm{Y}_{\text{batch}}$ that is difficult to map back to a state sequence.

The goal of the third step \textit{inverse problem solving} is the inverse mapping from $\bm{Y}_{\text{batch}}$ to $\bm{X}_{\text{batch}}$ using Bayesian estimation. We utilize the estimation framework of expectation propagation (EP) to approximate the posterior distribution of states with a Gaussian probability density function (pdf) with a mean $\widehat{\bm{X}}_{\textrm{batch}}$ and a variance $\widehat{\sigma}^2_{x, \textrm{batch}}$. The motivation for using EP over other Bayesian methods is its ability to handle non-Gaussianity priors on the states and provide a global approximation of the posterior through moment matching. Finally, in the \textit{uncertainty quantification} stage, we define the error uncertainty as a function of the variance $\widehat{\sigma}^2_{x, \textrm{batch}}$, which acts as an indicator of predictive confidence. The estimation of second-order moments is the cornerstone of this error quantification framework. By quantifying the spread of the estimated posterior through the EP, we move from the conventional point-estimate into a probabilistic approach, which is essential for error quantification and risk-aware control in general.

\subsection{Error quantification via Bayesian posterior variance}
While conventional data-driven methods typically provide point-estimates of the state trajectory, they offer no intrinsic measure of reliability. By transitioning to a probabilistic approach for the inverse problem solving (cf. step 3 in Fig. \ref{fig:inverse-problem-framework}), we define a rigorous indicator of predictive confidence through the spread of the estimated posterior distribution within the expectation propagation framework.

\subsubsection{The expectation propagation framework}
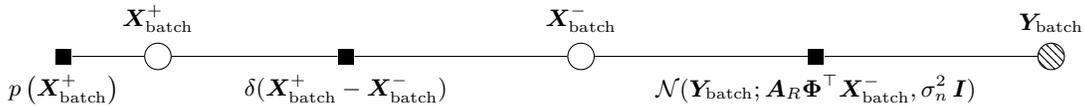
\begin{figure}[!h]
    \centering
    \begin{tikzpicture}[scale = 1.25]
    \newcommand{\vertex}{\node[vertex]}
    \tikzset{vertex/.style = {circle, draw, inner sep = 0pt, minimum size = 10pt}}
    \vertex[label = $\bm{X}_{\textrm{batch}}^{+}$](v1) at (-4.5,0) {};
    \vertex[label = $\bm{X}_{\textrm{batch}}^{-}$](v2) at (0,0) {};
    \vertex[label = {$\bm{Y}_{\textrm{batch}}$}, pattern = north west lines](z) at (5,0) {};
    \tikzset{vertex/.style = {rectangle, fill = black, inner sep = 0pt, minimum size = 6pt}}
    \vertex[label = below: $p\left(\bm{X}_{\textrm{batch}}^{+}\right)$](pv) at (-5.5,0) {};
    \vertex[label = below: $\delta(\bm{X}_{\textrm{batch}}^{+} - \bm{X}_{\textrm{batch}}^{-})$](delta) at (-2.5,0) {};
    \vertex[label = below: {$\mathcal{N}(\bm{Y}_{\textrm{batch}}; \bm{A}_R  \boldsymbol{\Phi}^\top\bm{X}_{\textrm{batch}}^{-}, \sigma^2_{n}\,\bm{I})$}](pw) at (2.5,0) {};
    \draw (pv)--(v1);
    \draw (delta)--(v1);
    \draw (delta)--(v2);
    \draw (pw)--(v2);
    \draw (pw)--(z);
    
\end{tikzpicture}
    \caption{Factor graph of the joint density factorization used in VAMP for error quantification.}
    \label{fig:VAMP-factor-graph-mmv}
\end{figure}

The objective of the inverse problem solving step is to resolve the inverse mapping from the predicted measurements $\bm{Y}_{\textrm{batch}}$ back to the underlying state sequence $\bm{X}_{\textrm{batch}}$. We achieve this by employing the Bayesian inference framework, namely expectation propagation described in Section \ref{sec:EP}. Fig. \ref{fig:VAMP-factor-graph-mmv} depicts the factor graph associated to the inverse problem in  (\ref{eq:Y=AX}) where the sensing matrix is $\bm{A}_R  \,\boldsymbol{\Phi}^\top$ and allowing an additive Gaussian noise $\bm{n}$ with variance $\sigma^2_n$.

Following the iterative message-passing process in Section \ref{sec:EP}, VAMP approximates the true posterior distribution of the states $\bm{X}_{\textrm{batch}}$ given the next states $\bm{Y}_{\textrm{batch}}$ as
\begin{equation}\label{eq:post-var}    p(\bm{X}_{\textrm{batch}}|\bm{Y}_{\textrm{batch}}) \approx \mathcal{N}(\bm{X}_{\textrm{batch}}; \widehat{\bm{X}}_{\textrm{batch}}, \widehat{\sigma}^2_{\bm{x},\text{batch}}),
\end{equation}
where $\widehat{\bm{X}}_{\textrm{batch}}$ is the posterior mean estimate of the states and $\widehat{\sigma}^2_{\bm{x},\text{batch}}$ its associated posterior variance.

\subsubsection{Variance as a metric for predictive confidence}

In the fourth and final \textit{uncertainty quantification} stage in Fig. \ref{fig:inverse-problem-framework}, the obtained posterior variance $\hat{\sigma}^2_{x,\textrm{batch}}$ in (\ref{eq:post-var}) is used as the primary metric for defining error uncertainty. This variance acts as a direct indicator of the ``recoverability'' of the state batch. Mathematically, if the data-driven model remains consistent with the true underlying state evolution, the inverse problem remains well-conditioned, resulting in a low posterior variance. In contrast, high variance values will highlight regions where the prediction of the data-driven model has entered a regime of high uncertainty, often due to a linearization failure in the high-curvature segments of the state manifold. This probabilistic characterization is essential for risk-aware control, as it allows the system to distinguish between high- and low-fidelity predictions where the deterministic output can no longer be trusted.

\section{Results} \label{ressec}
We present simulation results for three dynamical systems: two theoretical models namely the conductance-based neural model \cite{erme98} and the Hopf normal form model \cite{wigg03}, and one experimental system of a supersonic flow past a Cylinder \cite{wils20acc}. For our error quantification framework, we set the iteration number of VAMP to 50 and the states' prior $p(\bm{X}_{\textrm{batch}})$ to a Bernoulli Gaussian one with sparsity rate $\rho=5\%$.

\subsection{Conductance-Based Neural Model}\label{subsec:neural-model}
As a preliminary example, we consider a conductance-based neural model neuron \cite{wang96} with the addition of an adaptation current \cite{erme98}
\begin{align} \label{wbmodel}
C \dot{V}  &= -g_{\rm Na} m_\infty^3 q (V -E_{\rm Na}) - g_{\rm K} n^4(V -E_K) - g_{\rm L}(V-E_{\rm L}) - i_w  + u(t) + i_b, \nonumber \\
\dot{q} &= 5 \left[ \alpha_p(V)(1-q) - \beta_p(V)q   \right], \nonumber \\
\dot{n} &= 5 \left[ \alpha_n(V)(1-n) - \beta_n(V)n \right], \nonumber \\
\dot{w} &= a(1.5/(1+\exp((b-V)/k))-w).  
\end{align}
Above, $V$ is represents the transmembrane voltage, $q$ and $n$ are gating variables, and $i_w = g_w w (V-E_K)$ is an adaptation current governed by the variable $w$.  A constant baseline current $i_b = 10 \mu {\rm A}/{\rm cm}^2$ is used and the time-varying input $u(t)$ represents a transmembrane current.  The membrane capacitance, $C$, is taken to be $1 \mu  {\rm F}/{\rm cm}^2$.  Additional supporting equations are:
\begin{align*}
m_\infty &= \alpha_m(V)/(\alpha_m(V) + \beta_m(V)),   \\
\beta_n(V) &= 0.125\exp(-(V+44)/80),  \\
\alpha_n(V) &= -0.01(V+34)/(\exp(-0.1(V+34))-1), \\
\beta_q(V) &= 1/(\exp(-0.1(V+28))+1), \\
\alpha_q(V) &= 0.07\exp(-(V+58)/20), \\
\beta_m(V) &= 4\exp(-(V+60)/18), \\
\alpha_m(V) &= -0.1(V+35)/(\exp(-0.1(V+35))-1). 
\end{align*}
Reversal potentials and conductance are $E_{\rm Na} = 55 {\rm m}V, E_{\rm K} = -90{\rm m}V, E_{\rm L} = -65 {\rm  m}V, \, g_{\rm Na}= 35 {\rm mS}/{\rm cm}^2, g_{\rm K} = 9 {\rm mS}/{\rm cm}^2,  g_{\rm L} = 0.1 {\rm mS}/{\rm cm}^2, g_w = 2 {\rm mS}/{\rm cm}^2$. Auxiliary parameters are $a = 0.02\;{\rm ms}^{-1}, b = -5 \; {\rm m}V$ and $k = 0.5 {\rm m}V$.  When $u = 0$, Equation \eqref{wbmodel} fires periodic action potentials with period 6.53 ms. The input $u(t)$ modulates the firing rate of the action potentials.  In this example, let the state be $\bm{x}(t) = \begin{bmatrix} V(t) & q(t)  & n(t)  & w(t)  \end{bmatrix}^T$ with output
\begin{equation}
    g(\bm{x}) = \begin{bmatrix} V  \\ q \end{bmatrix}.
\end{equation}
Here, a constant time step is taken to be $\Delta t = 0.025$ so that $\bm{x}^k = \bm{x}(k \Delta t)$ and $u^k = u(k \Delta t)$.  In the model identification algorithm described in Sections \eqref{kooptheory} and \eqref{nonlinest}, a delay embedding length $z = 10$ is used.  A nonlinear lifting function $f_{\rm lift}(g(\bm{x}^k)) = f_1(f_2(g(\bm{x}^k)))$ is used where $f_2(g(\bm{x}^k)) \in \mathbb{R}^{10}$ is comprised of radial basis functions with the the $j^{\rm th}$ term given by $||g(\bm{x}^k) - \bm{q}_j||_2$ where $\bm{q}_j \in \mathbb{R}^2$ is the center of each radial basis function.  In this example, the first element of each $\bm{q}_j$ (associated with the voltage) is chosen randomly from a uniform distribution taking values between -300 and 200 and the second element (associated with the gating variable) is chosen randomly from a uniform distribution taking values between 0 and 1.  The function $f_1(f_2(g(\bm{x}^k))) \in \mathbb{R}^{990}$ applies a further lifting taking all polynomial combinations of the elements of $(f_2(g(\bm{x}^k)))$ up to degree 4.  Model identification is implemented using 300 time units of simulated data with the applied input $u(t) = 6 \sin(2 \pi t/200 + 0.0003t^2)$.

\begin{figure}[h!]
    \centering
    \hspace*{-0.9cm}
    \subfloat[state dynamics vs. $k \,\Delta t$\label{fig:simcc-state-eovlution}]{{
    \includegraphics[scale=0.235]{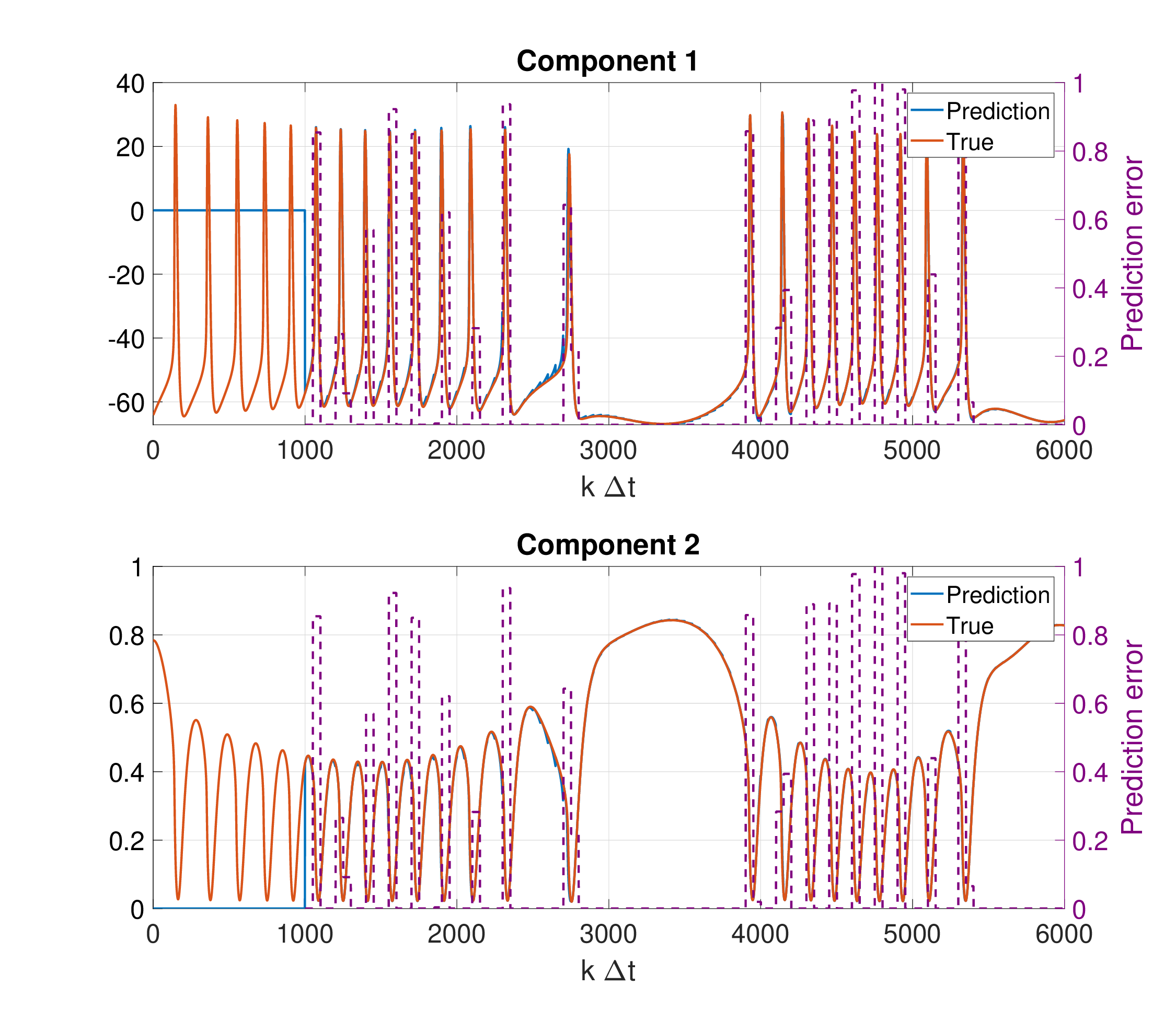} 
    }}%
    \subfloat[Uncertainty window vs. batch size\label{fig:simcc-batch-effect}]{{
    \includegraphics[scale=0.235]{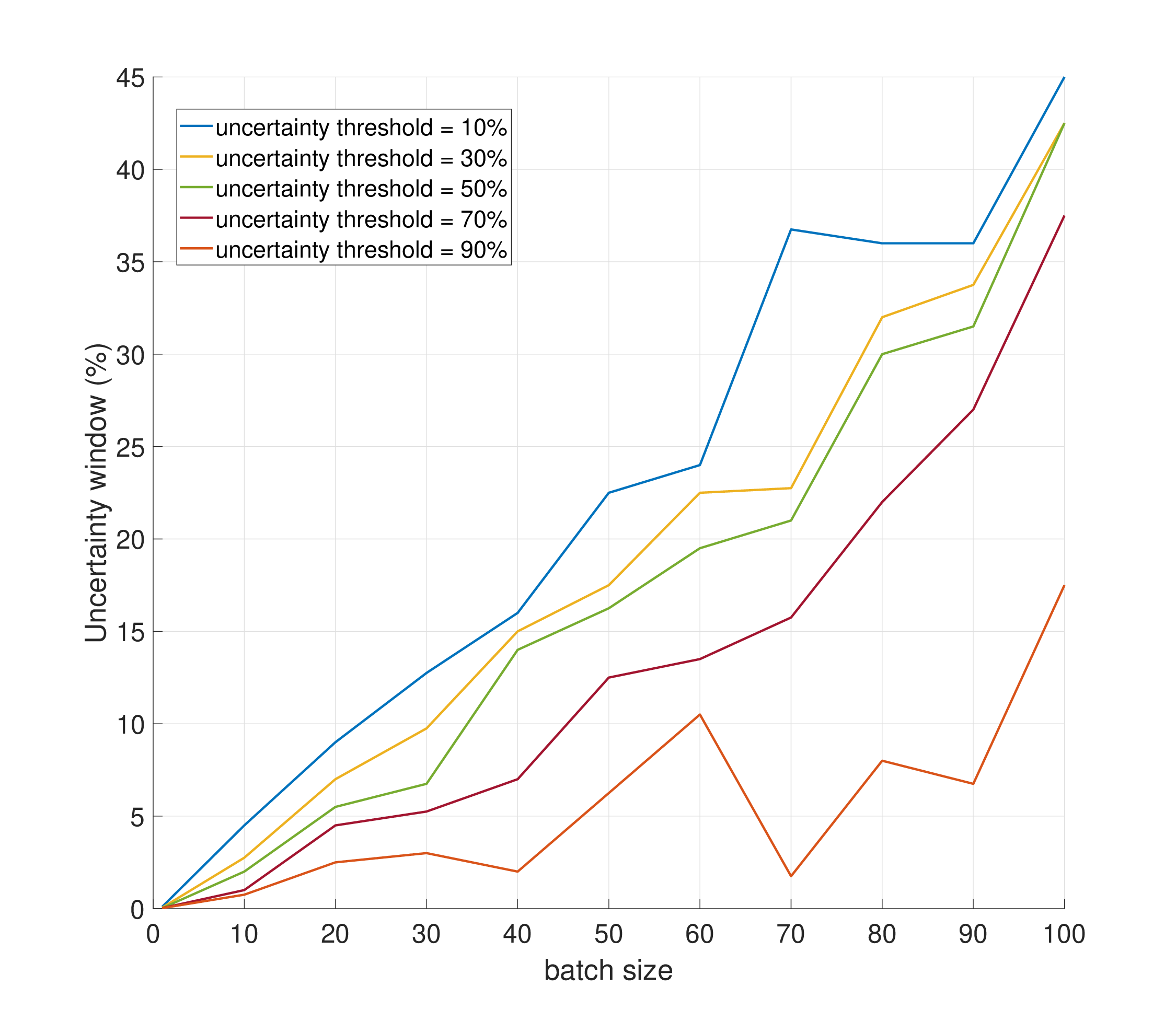} 
    }}%
    \caption{Uncertainty quantification of the predictive performance for the neural model.}
    \label{fig:simccmodel} 
\end{figure}

\begin{figure}[h!]
    \centering
    \includegraphics[scale=0.62]{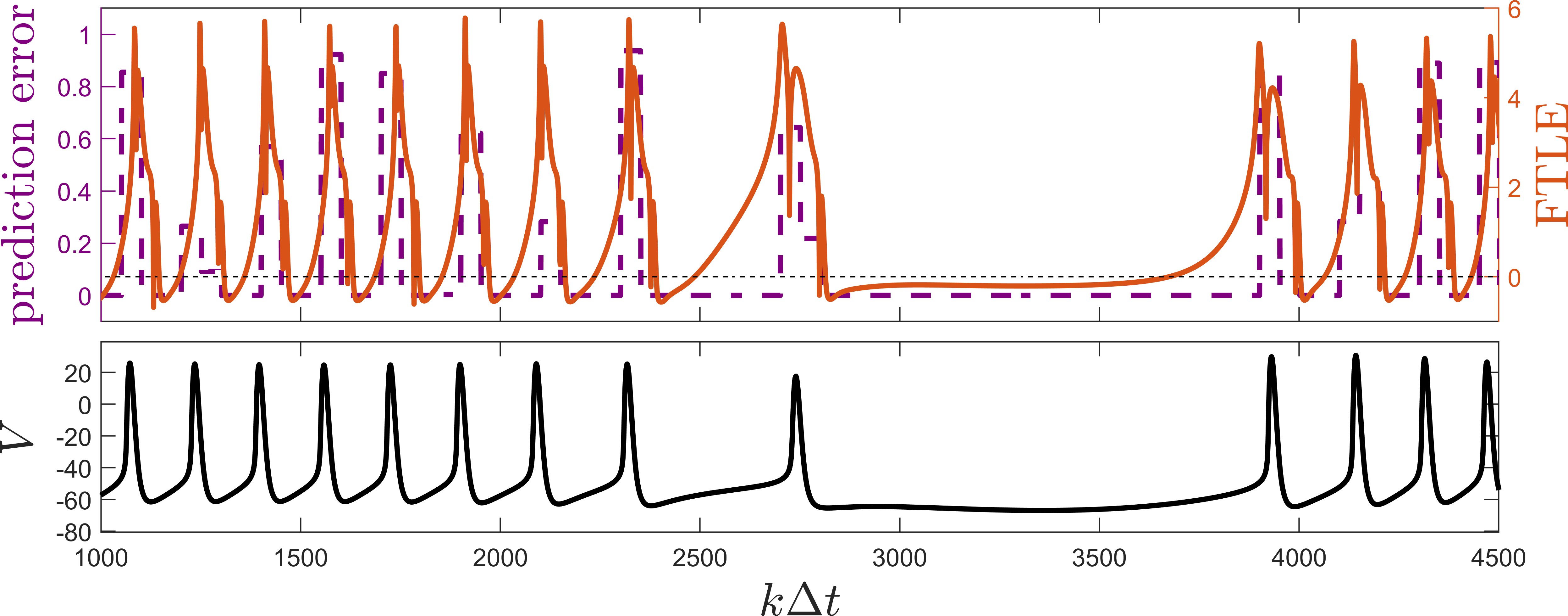}\vspace{-0.2cm}
    \caption{For the neural model \eqref{wbmodel}, the top panel shows the prediction error from Figure \ref{fig:simccmodel} superimposed on the FTLE, computed according to \eqref{leeq}.  The black dashed is shown for reference, corresponding to when the FTLE is zero and marking the boundary between amplification and attenuation of small perturbations to the state.  Over time, perturbations are amplified in regions where the FTLE is higher, explaining the larger prediction errors associated with these regions.  For reference, the bottom panel shows the transmembrane voltage highlighting that the prediction error and Lyapunov exponents are both large in the moments directly preceding an action potential.}
    \label{fig:leplot}
\end{figure}

In the prediction region (e.g., $k\Delta t > 1000$) in Fig. \ref{fig:simcc-state-eovlution}, it is seen that the data-driven model aligns closely with the true state components. However, as illustrated by the prediction error (in dashed purple) in Fig. \ref{fig:simcc-batch-effect}, the error exhibits significant volatility, often manifesting as high spikes when action potentials occur.  The behavior of the error can be explained by considering the finite time Lyapunov exponent computed over the preceding prediction window.  For the conductance-based neural model \eqref{wbmodel}, with $\bm{x} = [V \, q \, n \,w]^T$, for an infinitesimally small perturbation, $\bm{\delta x}$, the evolution over the prediction window, $\nu$, is given by
\begin{equation}
    \bm{\delta x}(t) = \bm{\Psi}(t,t-\nu) \bm{\delta x}(t-\nu) + O(|| \bm{\delta x} ||^2) ,
\end{equation}
where $\bm{\Psi}(t,t-\nu)$ is the state transition matrix associated with the linearization $\frac{d}{dt} \bm{\delta x} = \bm{J}(t) \bm{\delta x}$ with $\bm{J}(t)$ being Jacobian computed along the trajectory $\bm{x}(t)$.  Note that $\bm{\delta x}$ does not depend on the input in this example because it is simply appended to the voltage variable.  To leading order, the magnitude of this perturbation observed by the output is
\begin{equation}
    || \bm{\delta y}(t) ||^2 = \bm{\delta x}(t)^T \begin{bmatrix} \bm{e}_1 & \bm{e}_2 \end{bmatrix}^T  \bm{\Psi}^T(t,t-\nu) \bm{\Psi}(t,t-\nu) \begin{bmatrix} \bm{e}_1 & \bm{e}_2 \end{bmatrix} \bm{\delta x}(t),
\end{equation}
where $\bm{e}_1 = \begin{bmatrix} 1&0&0&0  \end{bmatrix}^T$, for instance.  A finite time Lyapunov exponent (FTLE) represents the maximum stretching seen from the perspective of the ouptput and is given by
\begin{equation} \label{leeq}
    \Lambda(t) = \frac{1}{\nu} \log( \sigma(t) ),
\end{equation}
where $\sigma$ is the maximum singular value of $ \bm{\Psi}(t,t-\nu) \begin{bmatrix} \bm{e}_1 & \bm{e}_2 \end{bmatrix} $.  Figure \ref{fig:leplot} shows the FTLE computed according to \eqref{leeq} superimposed over the prediction.  Regions where the FTLE is large correspond directly to regions for which the prediction error is also high.  The FTLE provides a local measurement of the sensitivity to initial conditions.  As such, uncertainties will be amplified in regions with higher FTLEs.

For any fixed uncertainty threshold, the prediction uncertainty window increases with the prediction horizon, which is directly related to the batch size as depicted in Fig. \ref{fig:simcc-batch-effect}. It is also observed how a low uncertainty threshold represents a risk-averse setting where the proposed method prioritizes error avoidance, effectively narrowing the model's predictive autonomy to ensure high fidelity. Conversely, higher thresholds represent a risk-tolerant regime that maximizes the duration of autonomous prediction by accepting higher levels of posterior variance associated to the inverse problem solution. It is also interesting to note that for higher uncertainty thresholds, the curve is non-monotonic indicating that the solution of the inverse problem's posterior variance is highly sensitive to the specific data points included in larger batches. Overall, the curves in Fig. \ref{fig:simcc-batch-effect} can be used to appropriately balance the tradeoff between prediction horizon and prediction uncertainty.

\subsection{Hopf Normal Form}
Here, we consider the Hopf normal form \cite{wigg03}
\begin{align} \label{hnf}
    \dot{x}_1 &= \sigma x_1(\mu - r^2) - x_2(1 + \rho(r^2-\mu)) + \sqrt{2 D} \eta(t) \nonumber \\
    \dot{x}_2 &= \sigma x_2 (\mu - r^2) + x_1(1 + \rho(r^2-\mu)), \nonumber \\
    g(x) &= x_1
\end{align}
Here, $\bm{x} = \begin{bmatrix} x_1 & x_2\end{bmatrix}^T$ where $x_1$ and $x_2$ are Cartesian coordinates, $y = x_1$ is the model output, $r^2 = x_1^2 + x_2^2$, and $\sqrt{2 D} \eta(t)$ is an independent and identically distributed zero-mean white noise process with intensity $D$.  This model undergoes a Hopf bifurcation at $\mu = 0$.  Here we take $\mu = 1$, $\rho = -0.1$, and $\sigma = 0.3$ so that this system has a stable limit cycle on the unit circle in the absence of noise.  In this example, the timestep is taken to be $\Delta t = 0.04$ with $\bm{x}^k = \bm{x}(k \Delta t)$.  In this example, there is no input and the system is driven by noise.  In the model identification algorithm described in Sections \ref{kooptheory} and \ref{nonlinest}, a delay embedding length $z = 9$ is used. The nonlinear lifting function $f_{\rm lift}(\bm{x}^k,\bm{h}^k) \in \mathbb{R}^{990}$ is used which takes polynomial combinations of the individual terms of $\bm{x}^k$ and $\bm{h}^k$ up to degree 4. 

\begin{figure}[h!]
    \centering
    \hspace*{-0.9cm}
    \subfloat[state dynamics vs. $k \,\Delta t$\label{fig:hopf-state-evolution}]{{
    \includegraphics[scale=0.23]{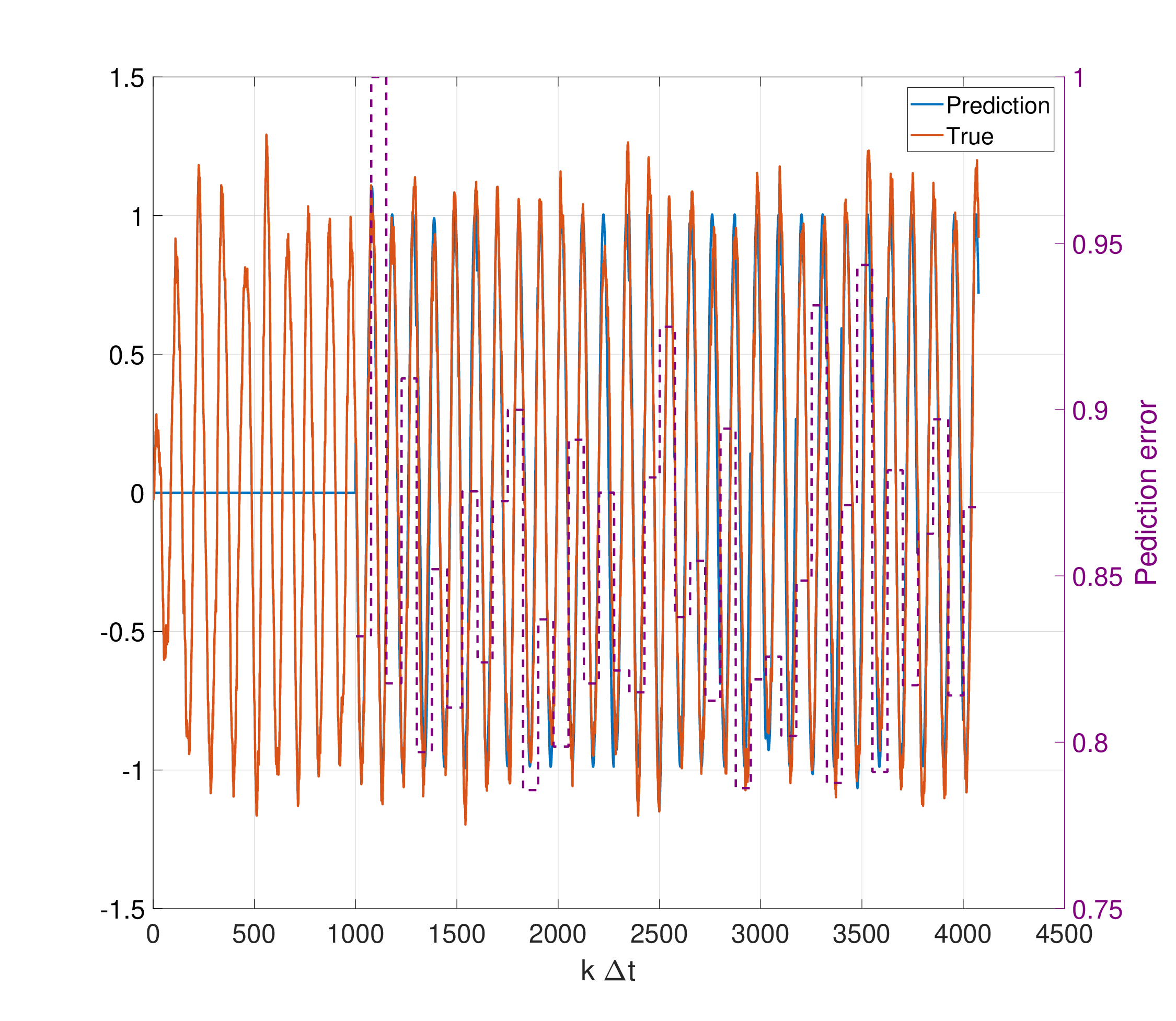} 
    }}%
    \subfloat[Uncertainty window vs. batch size\label{fig:hopf-batch-effect}]{{
    \includegraphics[scale=0.23]{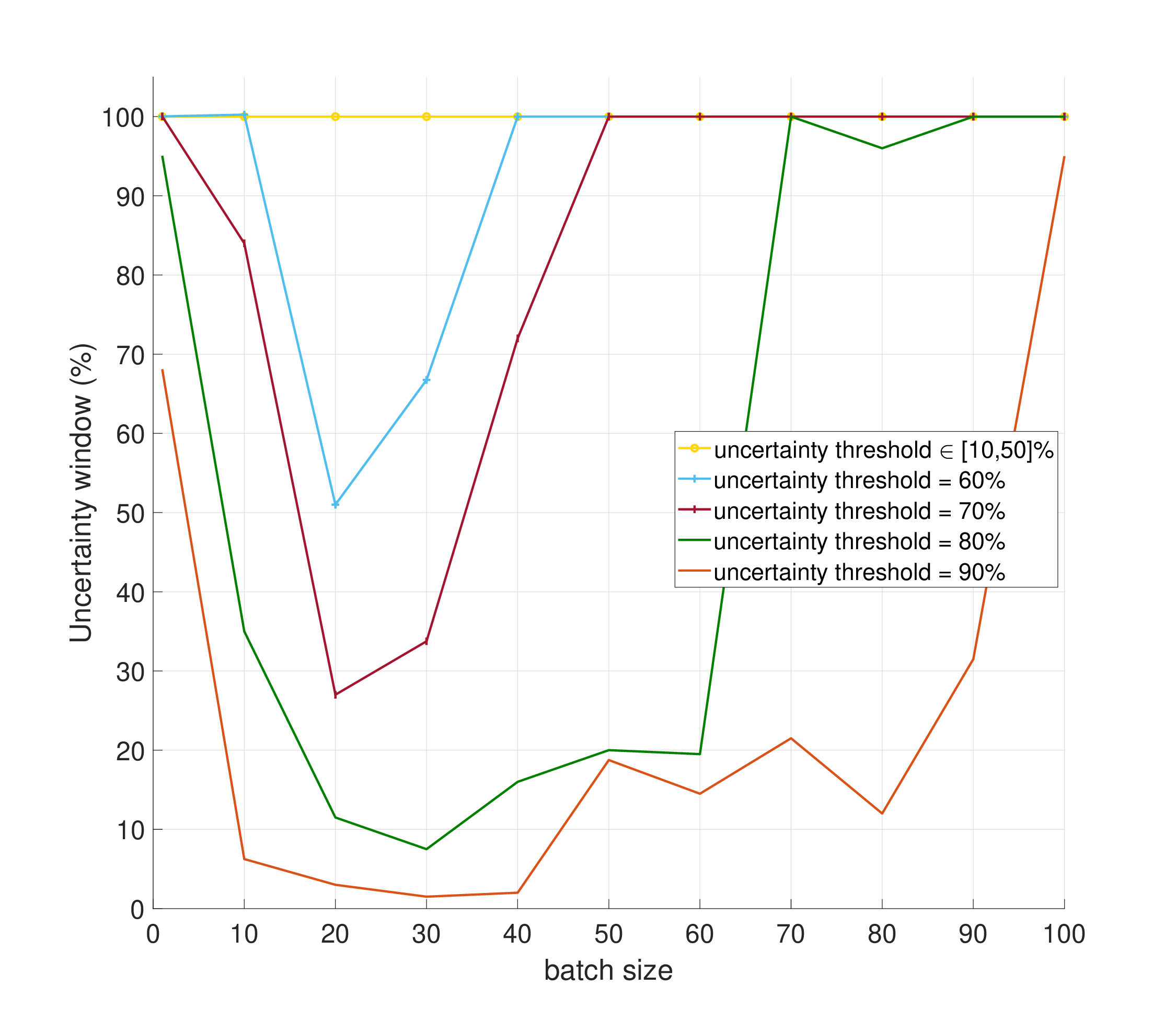} 
    }}%
    \caption{Uncertainty quantification of the predictive performance for the Hopf normal form \eqref{hnf}.}
    \label{fig:hopfmodel} 
\end{figure}

Unlike the neural model in Section \ref{subsec:neural-model}, the Hopf model does have a white noise process, which justifies why the prediction error (in dashed purple) in Fig. \ref{fig:hopf-state-evolution} maintains a high baseline  (between 0.8 and 1) and does not drops to lower values. There are visible pulses in the prediction error that align with the peaks of the true signal. This suggests that the error is state-dependent and occurs mainly around the extrema of the oscillations.

In Fig. \ref{fig:hopf-batch-effect}, for uncertainty thresholds between 10\% and 50\% (in yellow), the uncertainty window is flat at 100\% over the entire batch size range. This means the model is never trusted at these low levels of uncertainty thresholds. For higher thresholds, however, we see a U-shape where the uncertainty window decreases as the batch size moves from 0 to 30, before rising again. This suggest that very small batches provide insufficient data for the inverse problem to converge (i.e., high uncertainty due to lack of evidence), whereas very large batches suffer from drift (high uncertainty due to model cumulative prediction inaccuracy). It also means that there is a sweet spot for the batch size (around 20–30) where the Bayesian solver of the inverse problem is most confident.

\subsection{Supersonic Flow Past a Cylinder}
\begin{figure}[h!]
    \centering
    \hspace*{-0.9cm}
    \subfloat[Flow geometry \label{fig:lambdaedit}]{{
    \includegraphics[scale=0.235]{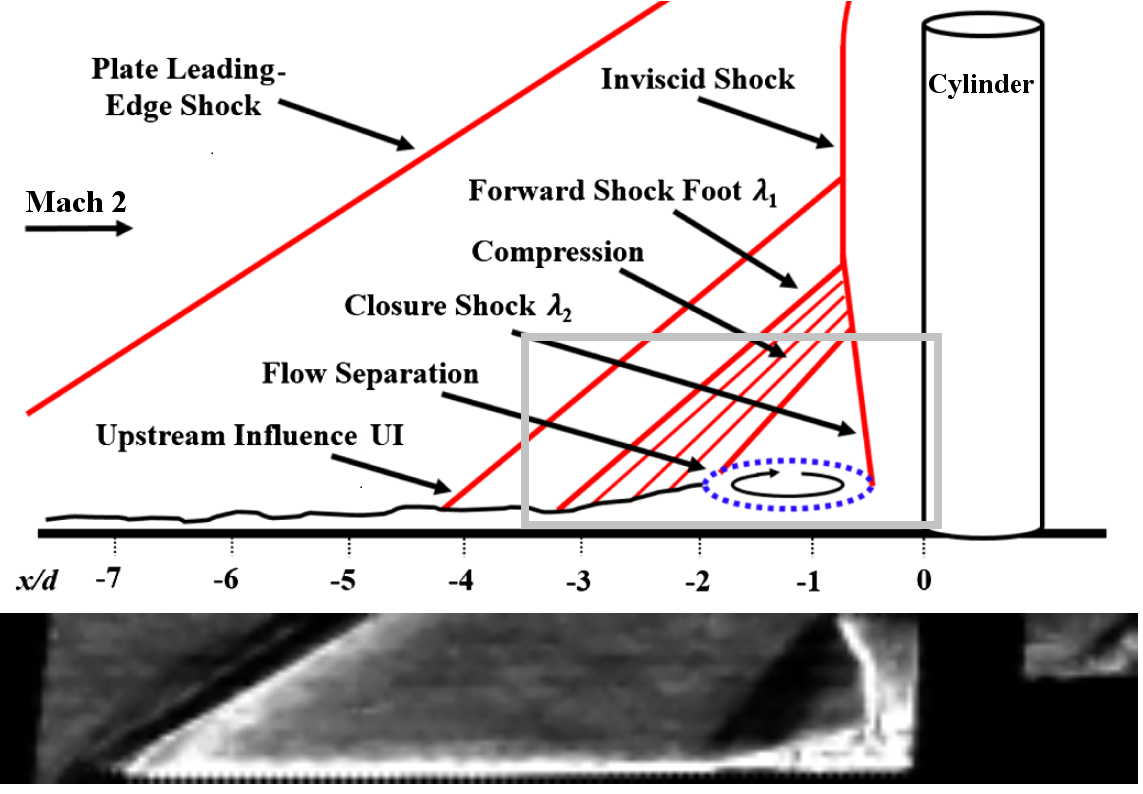} 
    }}%
    \subfloat[Dominant POD modes \label{fig:showmodes}]{{
    \includegraphics[scale=0.60]{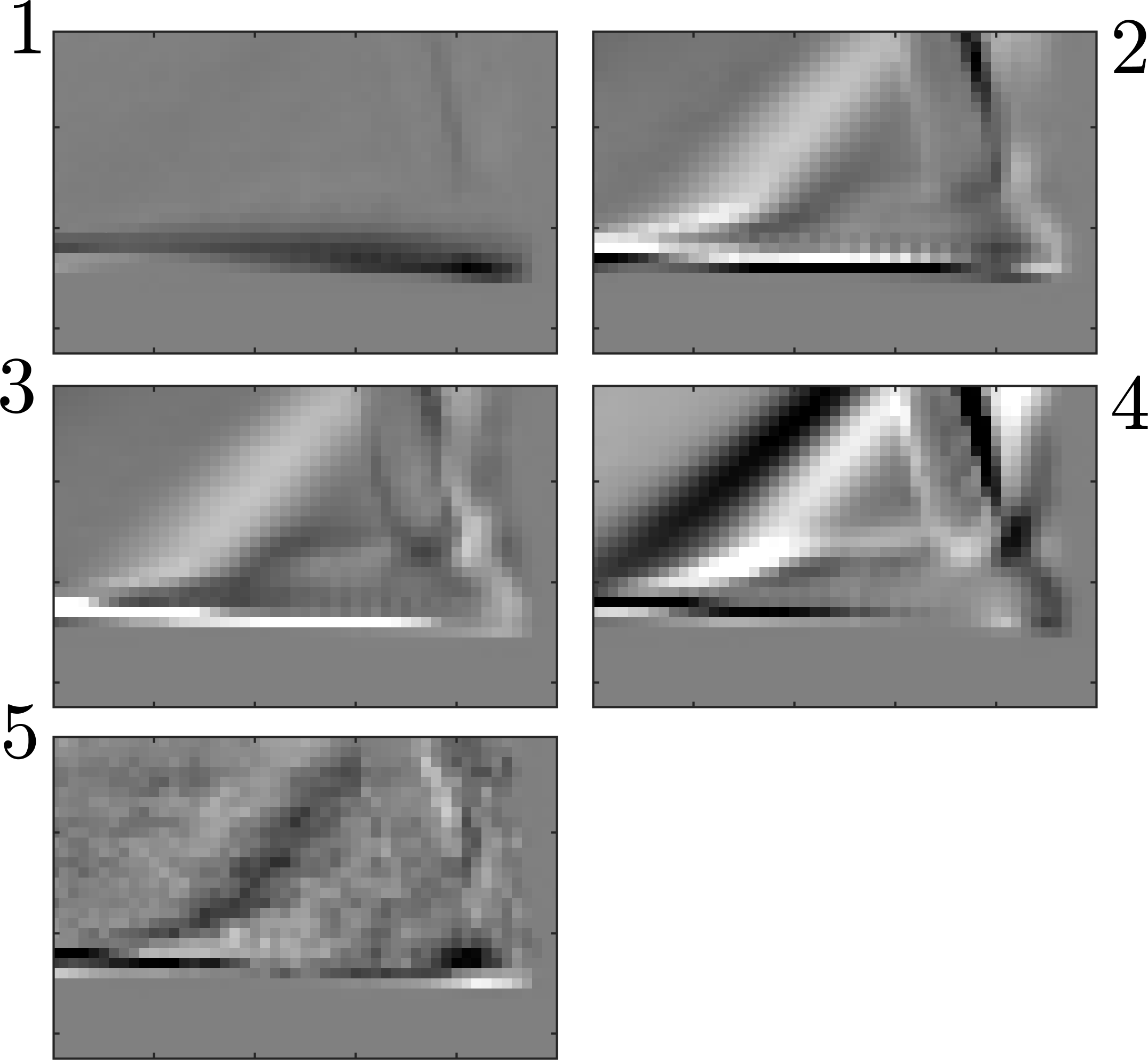} 
    }}%
    \caption{Panel (a) shows the flow geometry depicting shock-wave/transitional boundary layer-interaction in Mach 2 flow past the cylinder.  The black and white image is an example schlieren image.  Data analysis is performed on the region corresponding to the gray box.  The five most dominant POD modes are shown in panel(b).  Data is projected onto these modes and used in the data-driven model identification algorithm.}
    \label{fig:utsi-modes} 
\end{figure}

As a final illustration, we consider a set of experimental schlieren image data depicting cylinder-generated shock-wave/transitional boundary-layer interaction.  Relevant details of the experimental set up are given in \cite{wils20acc}.  Relevant features of the imaging data are depicted in Figure \ref{fig:utsi-modes} with panel (a) showing the flow geometry and an example schlieren image.  Here, pixel intensities (plotted in grayscale) correspond to the fluid density gradient.  The total data set contains 12,500 snapshots taken at 50 kHz with $\Delta t = 0.02$ ms between snapshots.  Each snapshot consists of 5,472 pixels; the combination of these pixel intensities is taken to be the state variable $\bm{x} \in \mathbb{R}^{5472}$.  Proper orthogonal decomposition is applied to the image data with Panel (b) showing the the five most dominant modes as gauged by the corresponding eigenvalues of the covariance matrix. These modes capture 48 percent of the total energy (i.e.~ taken as the sum of the 5 largest eigenvalues of the covariance matrix divided by the sum of the total eigenvalues).   The observable in this example is taken to be 
\begin{equation}
    g(\bm{x}^k) = \begin{bmatrix} \kappa_1^k & \dots & \kappa_5 \end{bmatrix}^T \in \mathbb{R}^5,
\end{equation}
where $\kappa_j^k$ corresponds to the $j^{\rm th}$ POD mode on the $k^{\rm th}$ snapshot.  We use a delay embedding of length $z = 25$ and a nonlinear lifting function $f_{\rm lift}(g(\bm{x}^k))$ comprised of all polynomial combinations of $g(\bm{x}^k)$ up to degree 4.

\begin{figure}[h!]
    \centering
    \hspace*{-0.9cm}
    \subfloat[state dynamics vs. $k \,\Delta t$\label{fig:utsi-state-evolution}]{{
    \includegraphics[scale=0.235]{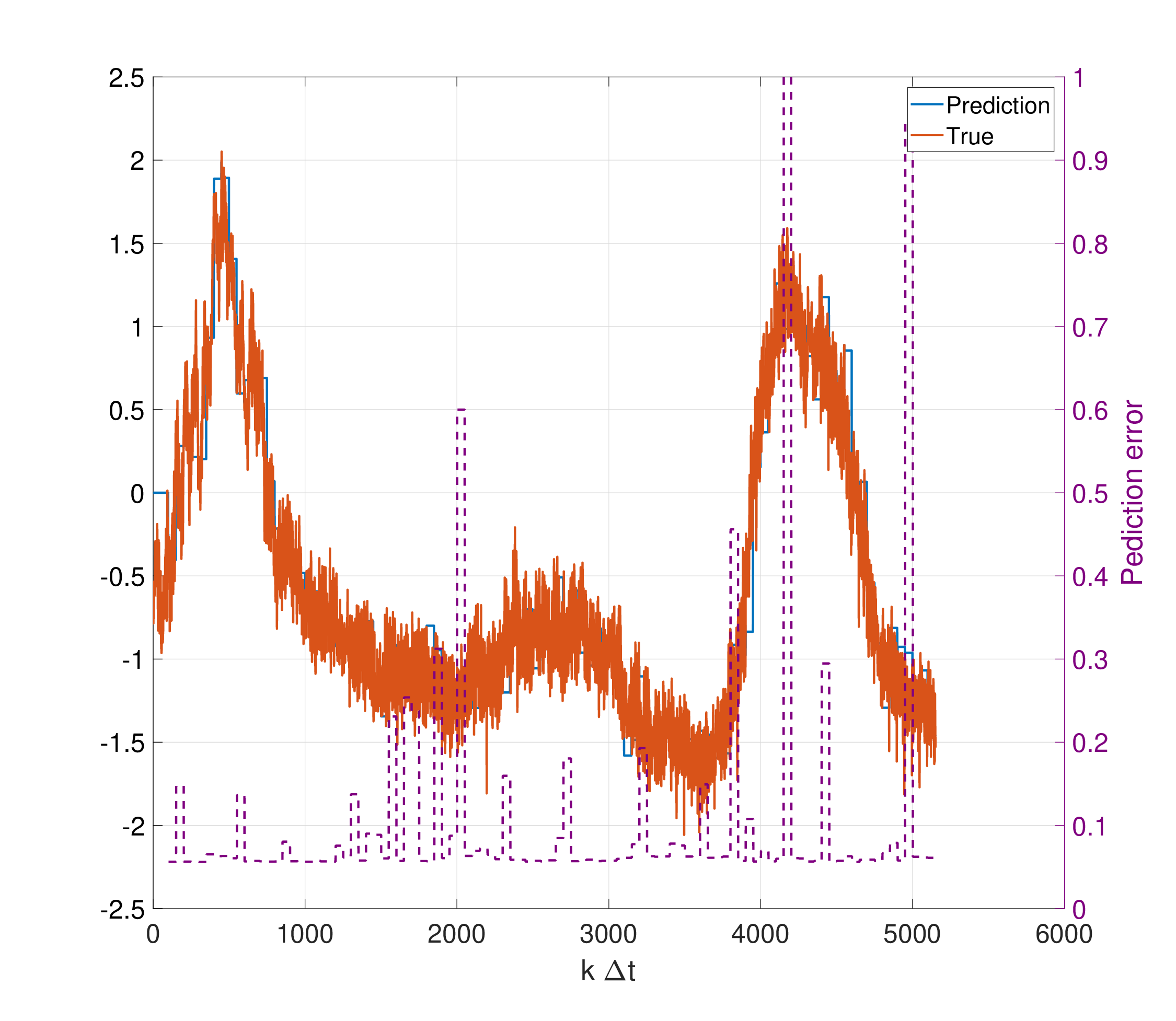} 
    }}%
    \subfloat[Uncertainty window vs. batch size\label{fig:utsi-batch-effect}]{{
    \includegraphics[scale=0.235]{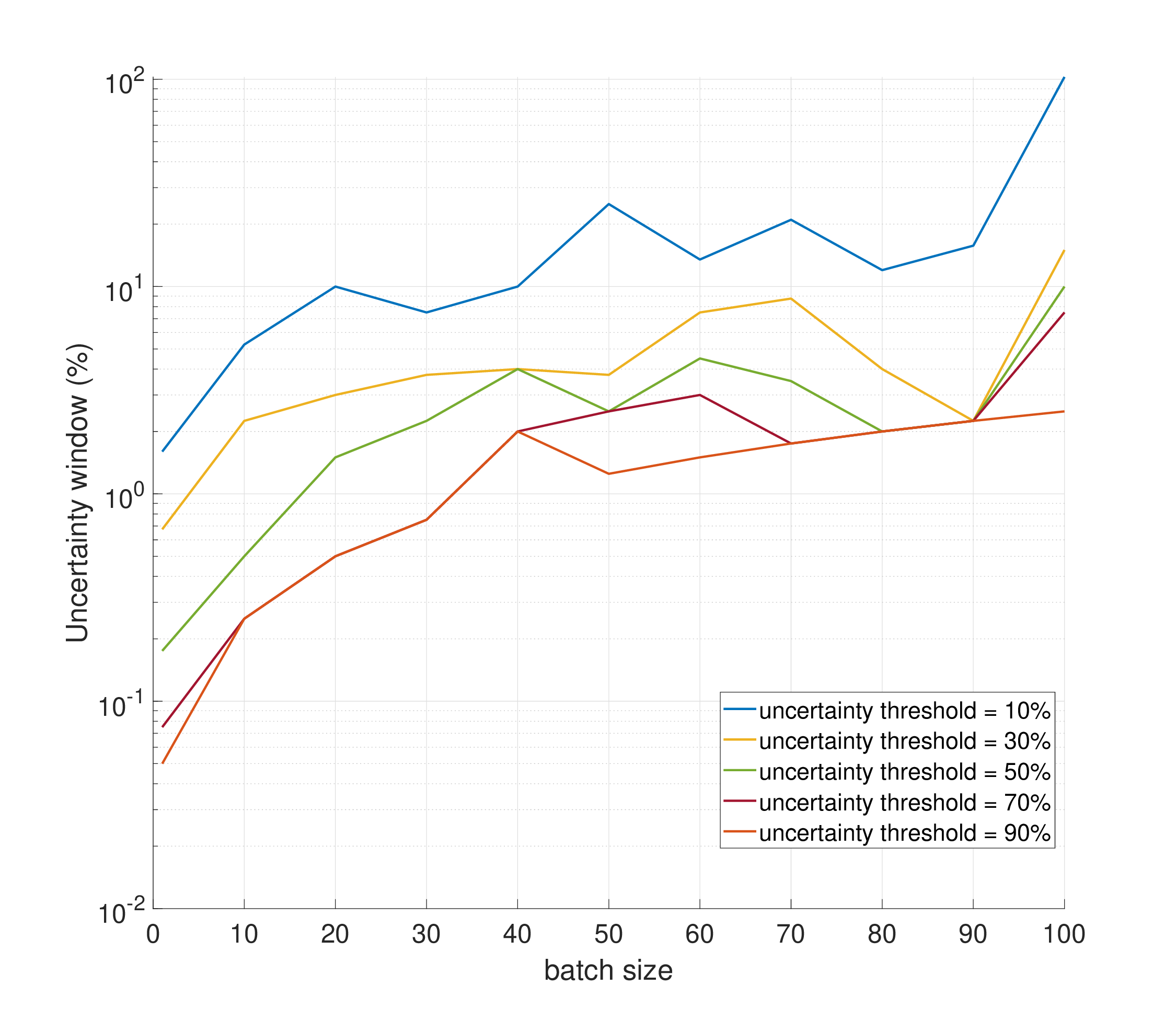} 
    }}%
    \caption{Uncertainty quantification of the predictive performance on the experimental schlieren imaging data of the supersonic flow past a cylinder.}
    \label{fig:utsi-model} 
\end{figure}

In the prediction region (e.g., $k\Delta t > 100$), the prediction of the data-driven model aligns closely with the true state components. However, as illustrated in Fig. \ref{fig:utsi-state-evolution}, the prediction error (in dashed purple) exhibits significant volatility, manifesting as high spikes during specific dynamical events. These error bursts typically correspond to moments when the signal enters multiple high-slope transitions.

Fig. \ref{fig:utsi-batch-effect} depicts the uncertainty window over batch sizes on a log scale. This indicates that the uncertainty window grows exponentially as the batch size increases from 0 to 100. It is observed that small batches ($<20$) and higher uncertainty thresholds yield very low uncertainty windows ($<1\%$). However, as the batch size grows, the probability of including states with significant inaccuracy (as the prediction of the data-driven model drifts) increases, causing the uncertainty window to jump by orders of magnitude.

The data-driven model prediction becomes unreliable very quickly as the horizon expands. Indeed, the performance degradation observed at extended prediction horizons (i.e., larger batch sizes) is primarily attributed to the inherent high noise level in the system, which makes the system less predictable within the inverse problem solving framework.

\FloatBarrier

\section{Conclusion} \label{concsec}
This work presents a new framework for error quantification in data-driven dynamical models  by recasting dynamical prediction as a Bayesian inverse problem. In doing so, our framework treats the predicted states from a data-driven model as pseudo-measurements and evaluates the predictive reliability without an external ground-truth oracle. Central to this approach is the expectation propagation (EP) algorithm for Gaussian message propagation which approximates the posterior distribution of states as a Gaussian message whose variance quantifies the prediction uncertainty. This framework relies on a linear inverse solver even though the the underlying dynamical model is non-linear. Future work could involve integrating the posterior variance into real-time adaptive control loops to trigger autonomous model re-synchronization or to dynamically tune system parameters based on the level of predictive confidence.  Additionally, more investigation on the fundamental relationship between finite time Lyapunov exponents and the prediction errors of the inferred models would be warranted.

This framework provides a principled bridge between data-driven model identification frameworks and Bayesian estimation. Its primary advantage is that it does not require an explicit analytical error model for the data-driven model and infers to uncertainty through the lens of recoverability from an estimation perspective. Furthermore, the use of EP ensures that the uncertainty estimates are more robust to outliers and non-linearities than standard linearized methods.

A limitation of this framework, however, is its computational demand. Solving an inverse problem for every prediction batch necessitates significant memory and processing power. Moreover, the framework is sensitive to the eigen spectrum of the sensing matrix $A$:~if the mapping between the states and the measurements is ill-posed, the resulting variance may reflect the weakness of the inverse setup rather than the failure of the data-driven model dynamics themselves.  Our error quantification framework is self-consistent: rather than comparing the model against an external oracle, we treat the model’s own forward prediction as a ``pseudo-measurement'' and attempt to invert it. This risky choice where the prediction serves as the ground truth for the inverse problem aligns with several established paradigms in signal processing and machine learning that utilize self-referential feedback to estimate reliability.

The most direct similarity of the proposed strategy is with temporal difference (TD) learning in reinforcement learning \cite{sutton1988learning}. In TD methods, the ground truth for the current state's value is not a final reward, but a prediction of the value of the next state, which is often called the TD target.  TD learning ``bootstraps'' by using one prediction to update another, our framework bootstraps by using the Koopman forward prediction $\bm{Y}_{\textrm{batch}}$ to quantify the error to reconstruct the state sequence $\bm{X}_{\text{batch}}$. Our error quantification framework shares the same strategy with autoencoders \cite{hinton2006reducing} pertaining to the concept of ``reconstruction loss''. An autoencoder compresses data and then attempts to reconstruct it. The ``error'' is measured by how much information was lost while compressing the input. Our framework treats the Koopman operator as the encoder for the forward state evolution and the Bayesian inverse solver as the decoder. Cycle consistency is used in computer vision to train models to translate images (e.g., from horse to zebra) by ensuring that translating the result back (e.g., from zebra to horse) returns the original image \cite{zhu2017unpaired}. Similarly, we are enforcing a consistency in time-series dynamics. If the forward mapping of the dynamic prediction and the backward mapping via inverse problem solving do not result in a high-fidelity recovery of the state, the model is deemed unreliable. While standard cycle consistency uses a point-wise error, we leverage expectation propagation to provide the estimation variance.

\section*{Acknowledgment} 
This material is based upon the work supported by the National Science Foundation (NSF) under Grant No.~CMMI-2024111 and the startup funding from the University of Tennessee.

\bibliography{myrefs}

\end{document}